\documentclass{amsart}

\usepackage{amsmath,amsthm}
\usepackage{amsfonts,amssymb}

\newtheorem{thm}{Theorem}[section]

\newtheorem{lem}[thm]{Lemma}
\newtheorem{prop}[thm]{Proposition}

\theoremstyle{remark}

 \def\cX{{\mathcal X}}
\newcommand{\PP}{{\mathsf {P}}}
\def\a{{\alpha}}
\def\b{{\beta}}

 \def\RR{{\mathbb R}}
 \def\ZZ{{\mathbb Z}}

 \def\supp{\operatorname{supp}}
 \def\sign{\operatorname{sign}}

\newcommand{\wh}{\widehat}

\begin{document}

\title
{Localized polynomial frames on the interval with Jacobi weights}

\author{Pencho Petrushev and Yuan Xu}
\address{Department of Mathematics\\University of South Carolina\\
Columbia, SC 29208.}
\email{pencho@math.sc.edu}
\address{Department of Mathematics\\ University of Oregon\\
    Eugene, Oregon 97403-1222.}\email{yuan@math.uoregon.edu}

\date{\today}
\thanks{The first author has been supported by NSF Grant DMS-0200665
and the second author by NSF Grant DMS-0201669}

\begin{abstract}
As is well known the kernel of the orthogonal projector onto the
polynomials of degree $n$ in $L^2(w_{\a,\b}, [-1, 1])$ with
$w_{\a,\b}(t) = (1-t)^\a(1+t)^\b$ can be written in terms of
Jacobi polynomials. It is shown that if the coefficients in this kernel
are smoothed out by sampling a $C^\infty$ function then the resulting function
has almost exponential (faster than any polynomial) rate of decay away
from the main diagonal.
This result is used for the construction of tight polynomial
frames for $L^2(w_{\a,\b})$ with elements having almost exponential
localization.
\end{abstract}

\maketitle

\section{Introduction}\label{Introduction}
\setcounter{equation}{0}

A basic technique in Harmonic analysis is
to represent functions or distributions as linear combinations of
functions of a particularly simple nature (building blocks),
which form bases or frames.
Meyer's wavelets \cite{Meyer} and the $\varphi$-transform
of Frazier and Jawerth \cite{F-J-W} provide such building blocks on $\RR^d$.
A distinctive feature of Meyer's wavelets and the elements of Frazier-Jawerth
is their almost exponential space localization and simple structure
on the frequency side.
This makes them an universal tool for decomposition of spaces of functions and
distributions on $\RR^d$.

Our primary goal in this article is to develop similar building blocks
for decomposition of weighted spaces on $[-1,1]$ with weight
$$
w_{\a,\b}(t) = (1-t)^\a(1+t)^\b, \qquad \a,\b > -1,
$$
and in particular, $L^p(w_{\a,\b})$, Hardy spaces $H^p(w_{\a,\b})$, Besov
spaces, and more general Triebel-Lizorkin spaces.
The structure of the weighted spaces on $[-1,1]$ is different and more
complicated than the structure of the spaces on $\RR$ due to the fact
that there are no simple translation or dilation operators in these spaces.
This creates a great deal of complications.

The Jacobi polynomials $\{P_n^{(\a,\b)} \}_{n=0}^\infty$
provide a basic vehicle for representation and analysis of functions or
distributions in the weighted spaces on $[-1,1]$ with weight
$w_{\a,\b}$.
We let $c_{\a,\b}$ denote the normalization constant of
$w_{\a,\b}$, i.e.
$c_{\a,\b}^{-1} := \int_{-1}^{1} w_{\a,\b}(t)dt$.
The Jacobi polynomials $\{P_n^{(\a,\b)}\}_{n=0}^\infty$ are orthogonal
with respect to $w_{\a,\b}$, namely \cite{Sz},
$$
  c_{\a,\b} \int_{-1}^1 P_n^{(\a,\b)}(t)
    P_m^{(\a,\b)}(t)w_{\a,\b}(t)dt  = \delta_{n, m} h_n^{(\a,\b)},
$$
where
$$
h_n^{(\a,\b)} = \frac{\Gamma(\a+\b+2)}{\Gamma(\a+1)\Gamma(\b+1)}
\frac{\Gamma(n+\a+1)\Gamma(n+\b+1)}{(2n+\a+\b+1)\Gamma(n+1)\Gamma(n+\a+\b+1)}.
$$
For $f \in L^1(w_{\a,\b})$ the Fourier expansion in Jacobi polynomials is
$$
  f(t) \sim \sum_{n=0}^\infty d_n(f) (h_n^{(\a,\b)})^{-1}
    P_n^{(\a,\b)}(t),
  \quad d_n(f) = c_{\a,\b}\int_{-1}^1 f(t) P_n^{(\a,\b)}(t)
            w_{\a,\b}(t)dt.
$$
The $n$th partial sum of the expansion can be written in terms of the
reproducing kernel $K_n^{(\a,\b)}(x,y)$ as
$$
  S_n(f;x) = \sum_{j=0}^n d_j(f) (h_j^{(\a,\b)})^{-1}
     P_j^{(\a,\b)}(x)
   = c_{\a,\b} \int_{-1}^1 f(t) K_n^{(\a,\b)}(x,t) w_{\a,\b}(t)dt,
$$
where the kernel is given by
\begin{equation}\label{eq:1.1}
K_n^{(\a,\b)}(x,y) = \sum_{j=0}^n
  \left(h_j^{(\a,\b)}\right)^{-1} P_j^{(\a,\b)}(x) P_j^{(\a,\b)}(y).
\end{equation}


One of our main results in this article asserts that any polynomial
in two variables of the form
\begin{equation}\label{def.L}
L_n^{\a,\b}(x,y)=\sum_{j=0}^\infty \wh a \left(\frac{j}{n}\right)
     \left(h_j^{(\a,\b)}\right)^{-1} P_j^{(\a,\b)}(x) P_j^{(\a,\b)}(y),
\end{equation}
where $\wh a \in C^\infty[0, \infty)$ with $\supp \, \wh a \subset [1/2, 2]$
(or more generally $\supp \wh a \subset [c, d]$, $c>0$) has almost exponential
localization around the main diagonal $y=x$ of $[-1, 1]^2$. To be more precise,
let
$$
w_{\a,\b}(n; x) := (1-x + n^{-2})^{\a+1/2} (1+x +
    n^{-2})^{\b+1/2},  \qquad -1 \le x \le 1,
$$
then for $\alpha, \beta > -1/2$ and any $k\ge 1$, there is a constant
$c_k > 0$ such that
\begin{equation}\label{L-loc}
|L_n^{\a,\b} (\cos\theta, \cos \phi)|
    \le c_k \frac{n}{\sqrt{w_{\a,\b} (n; \cos \theta)}
      \sqrt{w_{\a,\b} (n; \cos \phi)}
        (1+n|\theta - \phi|)^{k}}
\end{equation}
for $0 \le \theta, \phi \le \pi$. This result is a far reaching
extension of the recent discovery in \cite{NPW} that
$L_n^{\lambda,\lambda}(x,1)$ (with $\lambda$ a half integer)
has almost exponential localization around $x=1$,
which is utilized in \cite{NPW} for the construction of frames on
the $n$ dimensional sphere.
The above result shows that if one smooths out the coefficients
of the kernel in (\ref{eq:1.1}) by sampling a $C^\infty$ function
then the resulting function decays away from the main diagonal at
almost exponential (faster than any polynomial) rate.
This fact, which is well known for the trigonometric system,
has apparently been overlooked.
We believe that it will play an important role in various problems
for weighted spaces.

The polynomials $L_n^{\a,\b}$ give us a handy tool for
constructing tight frames in $L^2(w_{\a,\b})$
and other weighted spaces (with weight $w_{\a,\b}$) on $[-1, 1]$.
Our construction of frames utilizes a semi-discrete Calderon type decomposition
coupled with discretization using the Gaussian quadrature formula (see \S 3).
A~similar scheme is used in \cite{NPW} for the construction of frames
on the sphere.

Let us denote by
$\psi_\xi$, $\xi\in\cX$, the constructed frame elements (\S 3),
where $\cX=\cup_{j=0}^\infty \cX_j$ is a multilevel index set
consisting of the localization points (poles) of the $\psi_\xi$'s.
Then our frame is defined by
$$
\Psi:=\{\psi_\xi\}_{\xi\in\cX}.
$$
We show that every function $f \in
L^2(w_{\a,\b})$ has the representation
$$
f=\sum_{\xi \in \cX}\langle f, \psi_\xi \rangle \psi_\xi
\quad\mbox{and} \quad
\|f\|_{L^2(w_{\a,\b})}
=\Big(\sum_{\xi\in \cX} |\langle f, \psi_\xi\rangle|^2\Big)^{1/2},
$$
i.e. $\Psi$ is a tight frame in $L^2(w_{\a,\b})$.

A distinctive property of the $j$th level frame elements
$\psi_\xi$ ($\xi\in \cX_j$) is their localization:
\begin{equation}\label{psi-local}
|\psi_\xi(\cos \theta)| \le c_k
  \frac{2^{j/2}}{\sqrt{w_{\a,\b}(2^{j}; \cos \theta)}
     (1+2^j|\theta-\arccos \xi|)^k}
\qquad \forall \, k,
\end{equation}
which is analogous to the localization of Meyer's wavelets.
It is worthwhile to point out that the factor
$\sqrt{w_{\a,\b}(2^{j}; \cos \theta)}$ in the denominator in
(\ref{psi-local}) as well as the corresponding terms in (\ref{L-loc})
reflect the expected ``effect'' at the end points of $[-1, 1]$
and play a critical role in our development. Notice that this term is
present even in the case when $w_{\a,\b}(t)=1$ ($\a=\b=0$).

The superb localization of the frame elements $\{\psi_\xi\}$ prompted us
to term them {\bf needlets}.
This  
along with the semi-orthogonal structure of $\Psi$
and the increasing number of vanishing moments of the $\psi_\xi$'s
enables one to utilize the needlet system $\Psi$
for decomposition of spaces other than $L^2(w_{\a,\b})$
such as Besov and Triebel-Lizorkin spaces.
We shall report on results of this kind in a follow-up paper.

We do not discuss here the literature on polynomial bases and frames
on the interval because to the best
of our knowledge the elements of the existing bases or frames
do not have the localization of the needlets.
It is an open problem to construct bases in $L^2(w_{\a,\b})$
with basis elements having localization similar to (\ref{psi-local}).

This article is organized as follows.
In \S 2, we establish the localization (\ref{L-loc}) of the polynomials
$L_n^{\a,\b}$ from (\ref{def.L}),
where in \S 2.1 we consider the the particular case $\phi=0$ ($y=1$) ,
while in \S 2.2 we prove (\ref{L-loc}) in general.
In \S 3, we construct our polynomial frames and establish their
main properties.

Throughout the article positive constants are denoted by $c, c_1,\dots$;
unless specified, their values may vary at every occurrence. The notation
$A\sim B$ means $c_1A\le B\le c_2A$, and $A := B$ or $B =:A$ stands for
``$A$ is by definition equal to $B$''.

\section{Localized polynomials in terms of Jacobi polynomials}
\label{local-polyn}
\setcounter{equation}{0}

In this section we establish the localization properties of $L^{\a,\b}_n(x, y)$
from (\ref{def.L}) depending on the smoothness of $\wh a$.
We first prove (\ref{L-loc}) when $\phi=0$  ($y=1$)
and then consider the general case.

\subsection{The localization of \boldmath $L^{\a,\b}_n(x, y)$ in the case $y=1$}


\begin{thm} \label{thm:2.9}
Let $\wh a \in C^k[0,\infty)$ with $k\ge 1$ and $\supp \wh a \subset [1/2,2]$.
Assume that $\a \ge \b > -1/2$ and
define
\begin{equation}\label{def.Ln}
L_n^{\a,\b}(\cos \theta) = \sum_{j=0}^\infty \wh a\left(\frac{j}{n}\right)
  \left(h_j^{(\a,\b)}\right)^{-1} P_j^{(\a,\b)}(1)P_j^{(\a,\b)}(\cos \theta).
\end{equation}
Then there exists a constant $c_k>0$ depending only on $k$, $\a$, $\b$,
and $\wh a$ such that
\begin{equation}\label{est.Ln}
|L_n^{\a,\b}(\cos \theta)| \le c_k \frac{n^{2\a+2}}{(1+n \theta)^{k+\a-\b}},
  \quad 0 \le \theta \le \pi.
\end{equation}
The dependence of $c_k$ on $\wh a$ is of the form
$c_k=c(\a,\b,k)\max_{1\le\nu\le k}\|\wh a^{(\nu)}\|_{L^1}$.
\end{thm}

\begin{proof}
Since
$P_n^{(\a,\b)}(1) = (\a+1)_n / n!= \Gamma(n+\a+1)/[\Gamma(\a+1)\Gamma(n+1)]$
\cite[(4.1.1)]{Sz}, it is easy to verify that
\begin{align} \label{eq:L-n}
&L_n^{\a,\b}(\cos \theta) =  \frac{\Gamma(\b+1)}{\Gamma(\a+\b+2)} \\
   & \quad \quad \times  \sum_{j=0}^\infty \wh a\left(\frac{j}{n}\right)
     \frac{(2j+\a+\b+1)\Gamma(j+\a+\b+1)}{\Gamma(j+\b+1)}
   P_j^{(\a,\b)}(\cos \theta). \notag
\end{align}
Notice that since $\supp \wh a \subset [1/2,2]$,
the sum in (\ref{def.Ln}) (and (\ref{eq:L-n})) is finite;
in fact it is over $n/2<j<2n$.

We first prove (\ref{est.Ln}) for $0 \le \theta \le \pi/n$.
Using the fact that \cite[Theorem 7.32.1]{Sz}
$
\|P_j^{(\a,\b)}\|_{L^\infty[-1, 1]} \sim c j^{\alpha}
$
and $\Gamma(j + a) /\Gamma(j+1) \sim j^{a-1}$ as $j \to \infty$,
it follows that
\begin{align*}
|L_n^{\a,\b}(\cos\theta)|
     \le c \sum_{j=0}^{2n} j^{2 \a +1} \le c n^{2 \a+2},
\end{align*}
which yields (\ref{est.Ln}).

Let $\pi /n \le \theta \le \pi$. We will use the identity \cite{G}:
\begin{align*}
 & \frac{P_m^{(\a,\b)}(\cos \theta)}{P_m^{(\a,\b)}(1)} =
    d_{\a,\b} (1-\cos \theta)^{-\alpha}
 \int_0^\theta \cos \,[m+(\a + \b +1)/2] \phi  \\
   & \qquad\quad \times
  \frac{(\cos \phi - \cos \theta)^{\a-1/2}}{(1+\cos \phi)^{(\a+\b)/2}}
    {}_2F_1 \left( \frac{\a+\b}{2}, \frac{\a-\b}{2}; \a+\frac{1}{2};
     \frac{\cos \theta - \cos \phi}{1-\cos \theta}\right) d\phi,
\end{align*}
for $0 < \theta < \pi$ and $\alpha > -1/2$, where ${}_2F_1$ denotes
the hypergeometric function and $d_{\a,\b} = {2^{(\a+\b+1)/2}\Gamma(\a+1)}
/(\sqrt{\pi} \Gamma(\a+1/2))$.
Recall the identity \cite[(4.1.3)]{Sz}
$P_m^{(\a,\b)}(x) =  (-1)^mP_m^{(\b,\a)}(-x)$.
Now changing first the variables
$\phi \mapsto \pi - \phi$ in the integral, and then setting
$\theta \mapsto \pi - \theta$ and interchanging the places of
$\alpha$ and $\beta$, we obtain
\begin{align}\label{eq:jacobi}
 &\frac{P_m^{(\a,\b)}(\cos \theta)}{P_m^{(\b,\a)}(1)} =
    d_{\b, \a} (1+\cos \theta)^{-\b}
 \int_\theta^\pi \cos\,[m \phi -(\a + \b +1)(\pi-\phi)/2]  \\
   & \quad \times
  \frac{(\cos \theta - \cos \phi)^{\b-1/2}}{(1-\cos \phi)^{(\a+\b)/2}}
    {}_2F_1 \left( \frac{\a+\b}{2}, \frac{\b-\a}{2}; \b+\frac{1}{2};
     \frac{\cos \phi - \cos \theta}{1-\cos \theta}  \right) d\phi.
\notag
\end{align}
Since
$P_j^{(\b,\a)}(1) = \Gamma(j+\b+1)/[\Gamma(\b+1)\Gamma(j+1)]$,
we have
$$
P_j^{(\b,\a)}(1)
\frac{\Gamma(\b+1)}{\Gamma(\a+\b+2)}
 \frac{\Gamma(j+\a+\b+1)}{\Gamma(j+\b+1)}
= \frac{1}{\Gamma(\a+\b+2)}
 \frac{\Gamma(j+\a+\b+1)}{\Gamma(j+1)}.
$$
Using this and combining (\ref{eq:L-n}) with (\ref{eq:jacobi}), we get
\begin{align}\label{estim-Ln1}
&  L_n^{\a,\b}(\cos \theta) = d_{\a,\b}^* (1+\cos \theta)^{- \b}
    \int_\theta^\pi \left[\cos (\lambda \phi - \lambda \pi)A_n^{\cos}(\phi)
    - \sin (\lambda \phi - \lambda \pi)A_n^{\sin}(\phi) \right] \notag\\
  & \quad \times
  \frac{(\cos \theta - \cos \phi)^{\b-1/2}}{(1-\cos \phi)^{(\a+\b)/2}}
    {}_2F_1 \left( \frac{\a+\b}{2}, \frac{\b-\a}{2}; \b+\frac{1}{2};
     \frac{\cos \phi - \cos \theta}{1-\cos \theta}  \right) d\phi,
\end{align}
where $d_{\a,\b}^* := 2d_{\b,\a} / \Gamma(\a+\b+1)$ and
$$
A_n^{\cos}(\theta) := \sum_{j=0}^\infty \wh a \left(\frac{j}{n}\right)
\frac{(j+\lambda)\Gamma(j+2\lambda)}{\Gamma(j+1)} \cos j \theta,
$$
$$
A_n^{\sin}(\theta) := \sum_{j=0}^\infty \wh a \left(\frac{j}{n}\right)
\frac{(j+\lambda)\Gamma(j+2\lambda)}{\Gamma(j+1)} \sin j \theta
$$
with $\lambda :=(\alpha + \beta +1)/2$.

The idea of the proof is to derive (\ref{est.Ln})
from the analogous localization of the trigonometric polynomials
$A_n^{\cos}$ and $A_n^{\sin}$.
This in turn will follow by employing the fact that any
trigonometric polynomial
$
A_n(\theta) = \sum_{j=1}^\infty a_je^{ij\theta}
$
has an excellent localization around zero whenever the coefficients
${a_j}$ come from sampling of a smooth compactly supported function.

Define
$$
G(t):= \frac{(t+\lambda)\Gamma(t+2\lambda)}{\Gamma(t+1)}
\quad\mbox{if $t\ge0$.}
$$
Notice that $2\lambda=\a+\b+1 >0$, since $\a\ge \b>-1/2$.

With the next lemma we establish the smoothness and decay as $t\to\infty$
of $G$ and its derivatives.


\begin{lem} \label{lem:symptG}
The function $G$ is analytic on $[0, \infty)$ and for any $k\ge 0$
\begin{equation} \label{est.G}
|G^{(k)}(t)| \le c_k t^{2\lambda-k},
\quad t\ge 1,
\end{equation}
where $c_k$ depends only on $k$ and $\lambda$.
\end{lem}

\begin{proof}
We first show that for $b>a>0$, $\sigma:= a+1-b$, and $t > 0$,
\begin{equation} \label{G1}
\frac{d^k}{dt^k}\frac{\Gamma(t+a)}{\Gamma(t+b)}
=(-1)^kB_0^\sigma(a)\frac{\Gamma(b-a+k)}{\Gamma(b-a)}
\frac{1}{t^{b-a+k}}(1+\mathcal{O}(t^{-1}))
\quad\mbox{as $t\to\infty$,}
\end{equation}
where $B_n^\sigma(u)$ (here $n=0$) are the generalized Bernoulli polynomials
defined in \cite[(C.4.2)]{AAR}.
We will slightly modify the argument from \cite{AAR}.

The relation between Gamma and Beta functions gives
\begin{align*}
 \frac{\Gamma(t+a)}{\Gamma(t+b)} & = \frac{1}{\Gamma(b-a)}
    \int_0^1 u^{t+a -1}(1-u)^{b-a-1} du \\
  & = \frac{1}{\Gamma(b-a)}
    \int_0^\infty e^{-s (t+a)}(1- e^{-s})^{b-a-1} ds.
\end{align*}
It is readily seen that the last integral is an analytic function for $t\ge 0$
and it can be differentiated with respect to $t$ as follows
$$
\frac{d^k}{dt^k} \frac{\Gamma(t+a)}{\Gamma(t+b)}
  = \frac{(-1)^k}{\Gamma(b-a)}
    \int_0^\infty e^{-st} s^k e^{-sa}(1- e^{-s})^{b-a-1} ds.
$$
Now exactly as in \cite[C.4, p. 615]{AAR},
$$
\frac{s^ke^{-as}}{(1-e^{-s})^\sigma}
=\sum_{n=0}^\infty\frac{(-1)^n}{n!}B_n^\sigma(a)s^{n-\sigma+k},
\quad |t|<2\pi,
$$
and by Watson's lemma, we obtain the asymptotic expansion
$$
\frac{d^k}{dt^k}\frac{\Gamma(t+a)}{\Gamma(t+b)}
\sim \sum_{n=0}^\infty\frac{(-1)^{k+n}}{n!}B_n^\sigma(a)
\frac{\Gamma(b-a+k+n)}{\Gamma(b-a)}
\frac{1}{t^{b-a+k+n}}
\quad\mbox{as $t\to\infty$,}
$$
which implies (\ref{G1}).

Let $a := 2\lambda -  \lfloor 2\lambda \rfloor$, where
$\lfloor x \rfloor$ denotes the largest integer less than $x$.
Then $0 < a \le 1$.
Using that $\Gamma(1+x)=x\Gamma(x)$, we have
\begin{align*}
G(t)&=\frac{(t+1)(t+\lambda)\Gamma(t+ 2\lambda)}{\Gamma(t+2)}\\
&=(t+1)(t+\lambda)(t+2\lambda-1)\dots (t+a)\frac{\Gamma(t+a)}{\Gamma(t+2)}
=p(t)\frac{\Gamma(t+a)}{\Gamma(t+2)},
\end{align*}
where $p(t)$ is a polynomial of degree $\lfloor 2\lambda \rfloor+2$.
Now by (\ref{G1}),
$$
\left|\frac{d^\nu}{dt^\nu}\frac{\Gamma(t+a)}{\Gamma(t+2)}\right|
\le \frac{c}{t^{2-a+\nu}},
\quad t\ge 1, \quad\nu=0, 1, \dots,
$$
and hence
\begin{align*}
|G^{(k)}(t)| &\le \sum_{\nu=0}^k \binom{k}{\nu}
|p^{(\nu)}(t)|\Big|\frac{d^{k-\nu}}{dt^{k-\nu}}
\frac{\Gamma(t+a)}{\Gamma(t+2)}\Big|\\
&\le c\sum_{\nu=0}^k \binom{k}{\nu}
\frac{t^{\lfloor 2\lambda \rfloor+2-\nu}}{t^{2-a+k-\nu}}
\le ct^{2\lambda-k}.
\end{align*}
Thus (\ref{est.G}) is established.
\end{proof}

Now relying on the smoothness and localization properties of
$\wh a$ and $G$, we derive the desired
localization of $A_n^{\cos}$ and $A_n^{\sin}$.


\begin{lem} \label{lem:locA}
If $A_n=A_n^{\cos}$ or $A_n=A_n^{\sin}$, then
\begin{equation} \label{local.A}
|A_n(\theta)| \le c_k \frac{n^{2\lambda+1}}{(1+n|\theta|)^k},
\quad |\theta|\le \pi,
\end{equation}
where $c_k$ is of the form
$c_k=c(\a,\b,k)\max_{1\le\nu\le k}\|\wh a^{(\nu)}\|_{L^1}$
with $c(\a,\b,k)$ depending only on $\a$, $\b$, and $k$.
\end{lem}

\begin{proof}
Let $n\ge 2$ and set
$$
\Lambda_n(\theta):= \sum_{j\in \ZZ}
\wh a(|j|/n)s(j)G(|j|)e^{ij\theta},
$$
where $s(t)=1$ or $s(t)=\sign t$.
Evidently,
$\Lambda_n(\theta)=2A_n^{\cos}(\theta)$ if $s(t)=1$ and
$\Lambda_n(\theta)=2iA_n^{\sin}(\theta)$ if $s(t)=\sign t$.

Set
$\wh\Phi_n(t):=\wh a(|t|/n)s(t)G(|t|)$
and define
$
\Phi_n(t):= \frac{1}{2\pi}\int_\RR
\wh\Phi_n(\xi)e^{it\xi}d\xi.
$
Clearly, $\wh\Phi_n(t) \in C^k(\RR)$ and
$\supp \wh \Phi_n \subset [-2n,-n/2]\cup[n/2,2n]$.
Consequently,
$$
t^k\Phi_n(t)= \frac{i^k}{2\pi}\int_\RR
\frac{d^k}{d\xi^k}\wh\Phi_n(\xi)e^{it\xi}d\xi
$$
and using (\ref{est.G}),
\begin{align*}
|t^k\Phi_n(t)| &\le
\int_\RR\Big|\frac{d^k}{d\xi^k}[\wh a(\xi/n)G(\xi)]\Big|d\xi\\
  &\le \sum_{j=0}^k \binom{k}{j} \frac{1}{n^j}
     \int_{\RR} | \wh a^{(j)}(\xi/n)
       G^{(k-j)}(\xi)| d\xi \\
  & = 2\sum_{j=0}^k \binom{k}{j} \frac{n}{n^j}
     \int_{1/2}^2 |\wh a^{(j)}(u)
       G^{(k-j)}(n u)| du \\
  & \le cn^{2\lambda-k+1} \sum_{j=0}^k \binom{k}{j}
  \int_{1/2}^2 |\wh a^{(j)}(u)| u^{2 \lambda+j - k} du\\
  & \le c n^{2\lambda - k+1}\max_{0\le j\le k}\|\wh a^{(j)}\|_{L^1}
  =c_k n^{2\lambda - k+1}.
\end{align*}
Therefore,
$|\Phi_n(t)| \le c_k n^{2\lambda - k+1}|t|^{-k}$.
As above (with $k=0$) but easier, we obtain
$|\Phi_n(t)| \le c n^{2\lambda+1}$.
Combining these, we get
\begin{equation} \label{est.Phi}
|\Phi_n(t)| \le c_k\frac{n^{2\lambda+1}}{(1+n|t|)^k},
\quad t\in \RR.
\end{equation}

Applying the Poisson summation formula
$$
  2 \pi \sum_{j \in \ZZ} f(2 \pi j) = \sum_{j \in \ZZ} \wh f(j)
$$
to the function
$f(t):= \Phi_n(t+\theta)$
(then $\wh f(t) = \wh \Phi_n(t) e^{ it \theta}$),
we get
$$
 \Lambda_n(\theta) =
 \sum_{j \in \ZZ} \wh \Phi_n(j) e^{ij\theta} = 2 \pi \sum_{j \in \ZZ}
           \Phi_n(2 \pi j + \theta).
$$
Now using \eqref{est.Phi}, we obtain for $|\theta| \le \pi$,
\begin{align*}
|\Lambda_n(\theta)|
  & = 2 \pi \Big|\sum_{j \in \ZZ} \Phi_n (2 \pi j + \theta) \Big|
   \le c \sum_{j \in \ZZ} \frac{n^{2\lambda+1}}{ (1+ n |2 \pi j+\theta|)^k}\\
  & \le c \frac{n^{2\lambda+1}}{ (1+ n |\theta|)^k}
     + c\sum_{j = 1}^\infty \frac{n^{2\lambda+1}}{ (1+ n (2j-1)\pi)^k}\\
  & \le c \frac{n^{2\lambda+1}}{ (1+ n |\theta|)^k}
     + c\frac{n^{2\lambda+1}}{ (1+ n )^k}
        \sum_{j = 1}^\infty j^{-k} \\
  & \le c \frac{n^{2\lambda+1}}{ (1+ n |\theta|)^k}
\end{align*}
and estimate (\ref{local.A}) follows.
\end{proof}

We are now in a position to complete the proof of estimate (\ref{est.Ln})
in the case $\pi/n\le \theta\le \pi$.
Since
$0 \le (\cos \theta -\cos \phi)/(1-\cos \theta) \le 1$ for
$\theta \le \phi \le \pi$,
the absolute value of the ${}_2F_1$ term under the integral in 
(\ref{estim-Ln1}) is bounded by a constant. 
Hence, using (\ref{estim-Ln1}) and Lemma \ref{lem:locA},
we infer
\begin{equation}\label{eq:L-estimate}
|L_n^{\a,\b}(\cos \theta)| \le c (1+\cos \theta)^{- \b}
    \int_\theta^\pi \frac{n^{\a+\b+2}}{(1+ n \phi)^k}
    \frac{(\cos \theta - \cos \phi)^{\b-1/2}}{(1-\cos \phi)^{(\a+\b)/2}} d\phi.
\end{equation}
In the following we will use the identities $\cos \theta - \cos \phi =
 2 \sin \frac{\theta + \phi}{2} \sin \frac{\theta - \phi}{2}$,
$1- \cos \phi = 2 \sin^2 \frac{\phi}{2}$, and $1+ \cos \phi =
2 \cos^2 \frac{\phi}{2}$.

Consider first the case $\pi /n \le \theta \le \pi /2$.
Then for $\theta\le \phi\le \pi$,
$$
  \frac{2\sqrt{2}}{3\pi} (\theta+\phi) \le \sin \frac{\theta+\phi}{2}
   \le \frac{\theta+\phi}{2} \quad \hbox{and} \quad
 \frac{\phi-\theta}{\pi}  \le \sin \frac{\phi-\theta}{2}
     \le \frac{\phi-\theta}{2}.
$$
Using these in \eqref{eq:L-estimate}, we get
\begin{align*}
 \left| L_n^{\a,\b}(\cos \theta) \right|
   &  \le c n^{\a+\b+2}\int_\theta^\pi \frac{(\phi^2 - \theta^2)^{\b - 1/2}}
         {(1+n\phi)^k \phi^{\a+\b}} d\phi\\
   & \le c n^{\a+\b+2} \frac{\theta^{2\b}} {(n\theta)^k \theta^{\a+\b}}
       \int_1^{\frac{\pi}{\theta}} \frac{(u^2 - 1)^{\b-1/2}}{u^{k+\a+\b}} du\\
   & \le c \frac{n^{2\a+2}}{(n \theta)^{k +\a-\b}}
    \int_1^{\infty} \frac{(u^2 - 1)^{\b-1/2}}{u^{k+\a+\b}} du
      \le c \frac{n^{2\a+2}}{(1+n \theta)^{k+\a-\b}}.
\end{align*}

In the case $\pi/2 \le \theta \le \pi$, estimate \eqref{eq:L-estimate}
evidently gives
$$
 \left| L_n^{\a,\b}(\cos \theta) \right|
  \le c \frac{n^{\a+\b+2}}{(1+n\theta)^k} (1+\cos \theta)^{-\b}
    \int_\theta^\pi \frac{(\cos \theta - \cos \phi)^{\b -1/2}}
       {(1-\cos \phi)^{\a+\b}} d\phi.
$$
Let $\theta' = \pi-\theta$. Then $0 \le \theta' \le \pi/2$. Hence, for
$0 \le \phi \le \theta'$, $(\theta' \pm \phi)/\pi \le
\sin [(\theta' \pm \phi)/2] \le (\theta' \pm \phi)/2$.
Consequently, substituting $\phi \mapsto \pi - \phi$ gives
\begin{align*}
    \int_\theta^\pi \frac{(\cos \theta - \cos \phi)^{\b -1/2}}
       {(1-\cos \phi)^{\a+\b}} d\phi
& = \int_0^{\theta'}
    \frac{(\cos \phi - \cos \theta')^{\b -1/2}}
       {(1+\cos \phi)^{\a+\b}} d\phi \\
& \le c \int_0^{\theta'} ({\theta'}^2 - \phi^2)^{\b-1/2}d\phi
   =  c {\theta'}^{2\b}.
\end{align*}
Since $1+ \cos \theta = 1 - \cos \theta' = 2 \sin^2 (\theta'/2) \sim
(\theta')^2$, this shows that
$$
\left| L_n^{\a,\b}(\cos \theta) \right|
  \le c \frac{n^{\a+\b+2}}{(1+n\theta)^k} \le
     c \frac{n^{2\a+2}}{(1+n\theta)^{k + \a - \b}}
$$
as $\pi/2\le \theta \le \pi$.
The proof of Theorem~\ref{thm:2.9} is complete.
\end{proof}

\subsection{The localization of \boldmath $L^{\a,\b}_n(x, y)$
in the general case}\label{local-kernel}

Recall the definition of $L_n^{\a,\b}(x,y)$ from \eqref{def.L}:
$$
L_n^{\a,\b}(x,y):=\sum_{j=0}^\infty \wh a \left(\frac{j}{n}\right)
     \left(h_j^{(\a,\b)}\right)^{-1} P_j^{(\a,\b)}(x) P_j^{(\a,\b)}(y).
$$
In this subsection we estimate the localization
of $L_n^{\a,\b}(x, y)$ around the main diagonal $y=x$ of $[-1, 1]^2$,
which depends on the smoothness of $\wh a$.
To this end we will need the quantity
\begin{equation} \label{def.w}
w_{\a,\b}(n; x) := (1-x + n^{-2})^{\a+1/2} (1+x +
    n^{-2})^{\b+1/2},  \qquad -1 \le x \le 1.
\end{equation}
Notice that since
$\sin \tfrac{\theta}{2} \sim \sin \theta$ on $[0, 2\pi/3]$
and
$\cos \tfrac{\theta}{2} \sim \sin \theta$ on $[\pi/3, \pi]$,
then
\begin{equation} \label{w-sim1}
w_{\a,\b}(n;\cos \theta)\sim (\sin\theta + n^{-1})^{2 \a +1},
\quad 0\le \theta \le 2\pi/3,
\end{equation}
and
\begin{equation} \label{w-sim2}
w_{\a,\b}(n;\cos \theta)\sim (\sin\theta + n^{-1})^{2 \b +1},
\quad \pi/3\le \theta \le \pi.
\end{equation}


\begin{thm} \label{thm:3.1}
Let $\alpha, \beta > -1/2$ and let $\wh a \in C^k[0, \infty)$ with
$k \ge 2\a+2\b+3$ and $\supp \wh a \subset [1/2, 2]$.
Then there is a constant $c_k > 0$ depending only on
$k$, $\a$, $\b$, and $\wh a$ such that
for $0 \le \theta,\phi \le \pi$
\begin{equation} \label{Lbound1}
|L_n^{\a,\b} (\cos\theta, \cos \phi)|
    \le c_k \frac{n}{\sqrt{w_{\a,\b}(n; \cos \theta)}
     \sqrt{w_{\a,\b}(n; \cos \phi)} (1+n|\theta - \phi|)^{\sigma}},
\end{equation}
where $\sigma = k -2\a-2\b-3$.
Here the dependence of $c_k$ on $\wh a$ is of the form
$c_k=c(\a,\b,k)\max_{1\le\nu\le k}\|\wh a^{(\nu)}\|_{L^1}$.
\end{thm}

\begin{proof}
We need the product formula of Jacobi polynomials \cite{Koon}:
For $\a > \b > -1/2$,
\begin{align*}
\frac{P_n^{(\a,\b)}(x)  P_n^{(\a,\b)}(y)}{ P_n^{(\a,\b)}(1)}
 = c_{\a,\b} \int_{0}^\pi \int_0^1  P_n^{(\a,\b)} (t(x,y,r,\psi)) dm(r,\psi),
\end{align*}
where
$$
 t(x,y,r,\psi) = \tfrac{1}{2}(1+x)(1+y)+ \tfrac{1}{2} (1-x)(1-y) r^2 +
      r\sqrt{1-x^2}\sqrt{1-y^2}  \cos \psi -1,
$$
the integral is against
$$
  d m(r,\psi) = (1-r^2)^{\a-\b-1} r^{2\b+1} (\sin \psi)^{2\b} dr d\psi,
$$
and the constant $c_{\a,\b}$ is selected so that
$$
c_{\a,\b}\int_0^\pi\int_0^1 1 \,d m(r,\psi)=1.
$$
Once \eqref{Lbound1} is established for $\a >\b>-1/2$,
the case $\b > \a > -1/2$ will follow from the relation
$P_n^{(\a,\b)}(\cos \theta) = (-1)^n P_n^{(\b,\a)}(\cos (\pi -\theta))$.
In the case $\a = \b = \lambda -1/2$ we use the
product formula of Gegenbauer polynomials
\cite[Vol I, Sec. 3.15.1, (20)]{Edelyi}
$$
 \frac{C_j^\lambda(\cos \theta)C_j^\lambda(\cos \phi)}{C_j^\lambda(1)} =
    c_\lambda \int_{-1}^1 C_j^\lambda(\cos \theta \cos \phi + u\sin \theta
            \sin \phi) (1-u^2)^{\lambda -1} du,
$$
where $c_\lambda$ is the normalization constant of the weight function
$(1-u^2)^{\lambda -1}$. Since the case $\a = \b$ is much easier, we shall
only give the proof when $\a > \b$.

The product formula allows us to write $L_n^{\a,\b}(\cos \theta, \cos \phi)$
in terms of $L_n^{\a,\b}(t)$ defined in Theorem \ref{thm:2.9}, namely,
\begin{equation} \label{eq:Ln}
 L_n^{\a,\b}(x, y) = c_{\a,\b} \int_0^\pi \int_0^1
    L_n^{\a,\b}(t(x,y,r,\psi)) d m(r,\psi).
\end{equation}

For $t = \cos \theta$ with $0 \le \theta \le \pi$, we have $\theta  \sim
\sin \theta /2 \sim \sqrt{1-t}$,
and consequently the estimate of $|L_n^{\a,\b}(t)|$
from Theorem \ref{thm:2.9} can be rewritten as
\begin{equation} \label{est-Ln}
 \left| L_n^{\a,\b}(t) \right| \le
     c_k \frac{n^{2 \a +2}}{(1+n \sqrt{1-t})^{k+\a-\b}}, \quad -1 \le t \le 1.
\end{equation}
We will apply this estimate with $t= t(x,y,r,\psi)$. Let $x = \cos \theta$
and $y = \cos \phi$, $0 \le \theta, \phi \le \pi$. Then
\begin{align}\label{1-t}
  1-t(x,y,r,\psi) = &1- \cos\theta\cos\phi - \sin\theta\sin\phi \nonumber\\
     & +  \tfrac{1}{2} (1-\cos \theta)(1-\cos \phi) (1-r^2) +
      \sin \theta\sin\phi (1- r \cos \psi) \\
 = & 2 \sin^2 \frac{\theta -\phi}{2}
     + 2 \sin^2 \frac{\theta}{2}\sin^2 \frac{\phi}{2} (1-r^2) +
      \sin \theta\sin\phi (1- r \cos \psi),\nonumber
\end{align}
yielding
$$
1-t(x,y,r,\psi) \ge  2 \sin^2 \frac{\theta -\phi}{2} \sim (\theta - \phi)^2.
$$
Then by \eqref{eq:Ln}-\eqref{est-Ln},
we get
\begin{equation} \label{NoEndPt}
|L_n^{\a,\b}(\cos \theta,\cos \phi)| \le  \frac{ cn^{2 \a +2} }
   {(1+n|\theta - \phi|)^{k + \a - \b} }, \quad 0 \le \theta,\phi \le \pi.
\end{equation}

Our next step is to show that this estimate can be substantially
improved away from the end points of the interval.

\medskip\noindent
{\bf Case (a):} $0 \le \theta \le 2 \pi/3$ and $0 \le \phi \le \pi/2$.
By (\ref{1-t}), we have
$$
  1-t(x,y,r,\psi) \ge c (|\theta - \phi|^2 +
      \sin \theta \sin \phi (1-r \cos \psi)).
$$
Then by \eqref{eq:Ln}-\eqref{est-Ln},
\begin{align*}
|L_n^{\a,\b}(\cos \theta,\cos\phi)| \le & c n^{2\a+2}
\int_0^\pi \int_{0}^1 \frac{(1-r^2)^{\a-\b-1} r^{2\b+1} (\sin \psi)^{2\b}
     dr d\psi} {\left( 1+n[|\theta - \phi|^2 +
      \sin\theta \sin\phi
   (1-r \cos \psi)]^{1/2}\right)^\sigma},
\end{align*}
where $\sigma = k +\alpha-\beta$. Set $F(u) :=
1/( 1+n (|\theta - \phi|^2 +\sin\theta \sin\phi
(1-u))^{1/2})^\sigma$.
In the following we evaluate the double integral by
changing the variables several times:
\begin{align*}
\int_0^\pi \int_0^1 F(r \cos \psi) d m(r,\psi) & =
\int_0^1\int_{-1}^1 F(r s)(1-r^2)^{\a-\b-1} r^{2\b+1} (1-s^2)^{\b-1/2} dsdr\\
&=\int_0^1\int_{-r}^r F(u) r (1-r^2)^{\a-\b-1} (r^2-u^2)^{\b-1/2} du dr\\
&=\int_{-1}^1 F(u) \int_{|u|}^1 r (1-r^2)^{\a-\b-1} (r^2-u^2)^{\b-1/2}dr du\\
&= c \int_{-1}^1 F(u) (1-u^2)^{\a-1/2} du,
\end{align*}
where $c = (1/2) \int_0^1 v^{\b - 1/2} (1-v)^{\a-\b-1}dv$
(here $v=(r^2-u^2)/(1-u^2)$ is the last substitution).
This gives
\begin{align*}
|L_n^{\a,\b}(\cos \theta,\cos\phi)| \le & c n^{2\a+2}
 \int_{-1}^1 \frac{(1-u^2)^{\a-1/2} du}
  {( 1+n [|\theta - \phi|^2 +\sin\theta \sin\phi
(1-u)]^{1/2})^\sigma}.
\end{align*}
We now split up the last integral into two integrals:
one over $[-1,0]$ and the other over $[0,1]$.
For the first integral, $1- u \ge 1$ and hence
\begin{align*}
\int_{-1}^0 F(u) (1-u^2)^{\a-1/2} du & \le
\frac{c}{( 1+n [|\theta - \phi|^2 +\sin\theta \sin\phi]^{1/2})^\sigma}\\
& \le \frac{c}{(n^2\sin\theta \sin\phi)^{\a+1/2}
     (1+n|\theta-\phi|)^{k-\a-\b-1}}.
\end{align*}
To estimate the integral over $[0,1]$ we apply the substitution
$t=1-u$ and obtain
\begin{align*}
&  \int_0^1 F(u) (1-u^2)^{\a-1/2} du  \le
    c \int_0^1 \frac{(1-u)^{\a-1/2} du}{
(1+n[|\theta - \phi|^2 +\sin\theta\sin\phi
     (1-u)]^{1/2})^\sigma}\\
  & \qquad
    \le \frac{c}{(n^2 \sin\theta\sin\phi)^{\a+1/2}}
     \int_0^{n^2 \sin\theta\sin\phi}
      \frac{t^{\a-1/2}dt} {(1+[n^2(\theta-\phi)^2+t]^{1/2})^\sigma} \\
  & \qquad
    \le \frac{c}{(n^2 \sin\theta\sin\phi)^{\a+1/2}
      (1+n|\theta-\phi|)^{\sigma - 2\a -2} }
       \int_0^{\infty}   \frac{t^{\a-1/2} dt} {(1+t)^{\a+1}}\\
  & \qquad
    \le \frac{c}{(n^2 \sin\theta\sin\phi)^{\a+1/2}
      (1+n|\theta-\phi|)^{k - \a -\b -2} }.
\end{align*}
Putting these estimates together, we get
\begin{equation}\label{est-midle}
|L_n^{\a,\b}(\cos \theta,\cos\phi)| \le
  \frac{cn }{(\sin\theta\sin\phi)^{\a+1/2}
      (1+n|\theta-\phi|)^{k - \a -\b -2} }.
\end{equation}

Our last step here is to show that estimates
\eqref{NoEndPt}-(\ref{est-midle}) imply \eqref{Lbound1}.
Observe first that since $0 \le \theta, \phi \le 2\pi/3$,
we have by (\ref{w-sim1}),
\begin{align}\label{ww-sim}
w_{\a,\b}(n; \cos \theta)w_{\a,\b}(n; \cos \phi) &\sim
(\sin\theta+n^{-1})^{2\a+1}(\sin\phi+n^{-1})^{2\a+1}\\
&\sim (\theta+n^{-1})^{2\a+1}(\phi+n^{-1})^{2\a+1}.\notag
\end{align}
We need the following simple inequality
\begin{align}\label{simp-ineq}
(\theta+n^{-1})(\phi+n^{-1})
\le 3(\theta\phi+n^{-2})(1+n|\theta-\phi|),
\quad\theta, \phi\ge 0.
\end{align}
To prove this we assume that $\phi\ge\theta$ and
define $\gamma \ge1$ from $\phi=\gamma\theta$.
Suppose $\gamma\ge 3$. Then (\ref{simp-ineq}) will follow if
we show that $(\theta+\phi)n^{-1} \le 3n^{-2}n|\theta-\phi|$,
which is equivalent to $\gamma+1\le 3(\gamma-1)$.
But this holds since $\gamma\ge 3$.

Let now $1\le \gamma <3$. Then it suffices to show that
$(\theta+\phi)n^{-1} \le 2(\theta\phi+n^{-2})$.
In turn this inequality holds if
$4\theta n^{-1}\le 2(\theta^2+n^{-2})$,
which is apparently true.
Thus (\ref{simp-ineq}) is established.

It is now readily seen that \eqref{NoEndPt}-(\ref{est-midle})
yield \eqref{Lbound1} by using \eqref{ww-sim}-\eqref{simp-ineq}.


\medskip\noindent
{\bf Case (b):} $\pi/3 \le \theta \le \pi$ and $\pi/2 \le \phi \le \pi$.
In this case, we have by (\ref{1-t})
\begin{equation} \label{eq:1-t(b)}
 1-t(x,y,r,\psi) \ge c (|\theta - \phi|^2 + (1-r) +
       r(1-\cos \psi)\sin\theta\sin\phi)
\end{equation}
upon writing $1- r \cos \psi = 1- r + r (1-\cos \psi)$.
In particular,
$$
 1-t(x,y,r,\psi) \ge c (|\theta - \phi|^2 + (1-r)),
$$
which shows, using \eqref{eq:Ln}-(\ref{est-Ln}) , that
\begin{align*}
|L_n^{\a,\b}(\cos \theta,\cos\phi)| & \le  c n^{2\a+2}
  \int_0^1  \frac{(1-r^2)^{\a-\b-1} r^{2\b+1} dr}
{\left( 1+n[|\theta - \phi|^2 + (1-r)]^{1/2}\right)^\sigma},
\end{align*}
with $\sigma := k +\alpha-\beta$ as before.
Applying the substitution $n^2 (1-r) = s$, we get
\begin{align} \label{eq:NoEnd-b}
|L_n^{\a,\b}(\cos \theta,\cos\phi)| & \le  c n^{2\a+2}
   \int_0^1  \frac{(1-r)^{\a-\b-1} dr}
   {\left( 1+n (|\theta - \phi|^2 + (1-r))^{1/2} \right)^{\sigma}} \\
 & \le  c n^{2\b + 2}  \int_0^{n^2}  \frac{s^{\a-\b-1} ds}
   {\left( 1+ (n^2|\theta - \phi| + s )^{1/2} \right)^{\sigma}}\notag \\
 & \le  c \frac{ n^{2\b+2} } { (1+n|\theta - \phi|)^{\sigma - 2(\a-\b)-1/2 }}
    \int_0^{\infty}  \frac{s^{\a-\b-1} ds}
         {(1+ s)^{\a-\b+1/4}} \notag \\
 & \le  c \frac{ n^{2\b+2} } { (1+n|\theta - \phi|)^{k - \a +\b-1/2 }}. \notag
\end{align}
Furthermore, using \eqref{eq:Ln}-(\ref{est-Ln}) and \eqref{eq:1-t(b)},
we obtain
\begin{align*}
& |L_n^{\a,\b}(\cos \theta,\cos\phi)| \\
& \qquad \le  c n^{2\a+2}
\int_0^1 \int_0^\pi \frac{(1-r^2)^{\a-\b-1} r^{2\b+1} (\sin \psi)^{2\b}
     d\psi dr} {\left( 1+n[|\theta - \phi|^2 + (1-r) +
r(1-\cos \psi)\sin\theta\sin\phi]^{1/2}\right)^\sigma}.
\end{align*}
The inner integral can be estimated by changing the variables.
Set $A := |\theta - \phi|^2+1- r$. Then
\begin{align*}
 \int_0^{\pi/2} & = \int_{0}^1 \frac{(1-s)^{\b-1/2} ds}
{(1+n[A+r(1-s) \sin\theta\sin\phi]^{1/2})^\sigma}\\
& = \frac{1}{(n^2 r\sin\theta\sin\phi)^{\b + 1/2}}
   \int_0^{n^2 r \sin\theta\sin\phi}
     \frac{t^{\b -1/2} dt} { (1+ (n^2 A + t)^{1/2})^\sigma} \\
& \le \frac{1}{(n^2 r \sin\theta\sin\phi)^{\b + 1/2}
        (1+ n A^{1/2})^{ \sigma - 2\beta -2} }
   \int_0^{\infty}
     \frac{t^{\b -1/2} dt} { (1+ t)^{\beta+1} } \\
& \le \frac{c}{(n^2 r \sin\theta\sin\phi)^{\b + 1/2}
        (1+ n A^{1/2})^{ k+\a - 3\beta -2} },
\end{align*}
and since $1- s \ge 1$ for $-1 \le s \le 0$,
\begin{align*}
  \int_{\pi/2}^{\pi} & = \int_{-1}^0 \frac{(1-s)^{\b-1/2} ds}
     { (1+n (A+r(1-s) \sin\theta\sin\phi)^\sigma}
  \le \frac{c}{(1+n (A+r\sin\theta\sin\phi)^{1/2})^\sigma} \\
 &  \le \frac{c}{(n^2 r \sin\theta\sin\phi)^{\b + 1/2}
        (1+ n A^{1/2})^{ k+\a - 3\beta -2} }.
\end{align*}
Combining these, we get
\begin{align*}
 |L_n^{\a,\b}(\cos \theta,\cos\phi)|  \le c \frac{ n^{2\a+2}}
{(n^2\sin\theta\sin\phi)^{\b + 1/2}}
\int_0^1  \frac{(1-r^2)^{\a-\b-1} r^{\b+1/2}dr}
   {\left( 1+n (|\theta - \phi|^2 + (1-r))^{1/2} \right)^{\tau}},
\end{align*}
where $\tau := k +\a - 3\b -2$.
The last integral (with $\sigma$ in place of $\tau$) has already been
estimated in \eqref{eq:NoEnd-b}. Consequently,
$$
 |L_n^{\a,\b}(\cos \theta,\cos\phi)|  \le  \frac{cn^{2\b+2}}
{(n^2 \sin\theta\sin\phi)^{\b + 1/2}
  (1+n|\theta - \phi|)^{k - \a -\b -5/2}}.
$$
As in Case (a), this estimate along with (\ref{eq:NoEnd-b}) implies
(\ref{Lbound1}), using (\ref{w-sim2}) and \eqref{simp-ineq}.


\medskip\noindent
{\bf Case (c):} Either $0 \le \theta \le \pi/3$ and $\pi/2 \le \phi \le \pi$
or $2\pi/3 \le \theta \le \pi$ and $0 \le \phi \le \pi/2$.  In this case,
we have $|\theta - \phi| \ge \pi/6$, so that
$$
  1-t(x,y,r,\psi) \ge c |\theta - \phi|^2 \ge  c >0.
$$
Hence, by \eqref{eq:Ln}-(\ref{est-Ln}),
\begin{align*}
|L_n^{\a,\b}(\cos \theta,\cos\phi)|
\le \frac{c n^{2\a+2}}{(1+n)^{k + \a -\b}}
\le \frac{c n}{(1+n)^{k - \a -\b -1}},
\end{align*}
which yields (\ref{Lbound1}) since
$w_{\a,\b}(n;\cos\theta) w_{\a,\b}(n;\cos\phi) \le c$.
The~proof of Theorem~\ref{thm:3.1} is complete.
\end{proof}

The estimate of $|L_n^{\a,\b}(x,y)|$ from Theorem~\ref{thm:3.1}
allows us to control its $L^p$ integral.


\begin{prop} \label{prop:3.2}
Let $k \ge 3\a+3\b+5$. Then for $1 \le p < \infty$,
\begin{equation} \label{est-Lp-norm}
\int_{-1}^1 |L_n^{\a,\b} (x,y)|^p w_{\a,\b}(y) dy \le c
\frac{n^{p-1}}{(w_{\a,\b}(n;x))^{p-1}},
\quad -1\le x\le 1.
\end{equation}
\end{prop}

\begin{proof}
We will prove (\ref{est-Lp-norm}) only in the case $0 \le x \le 1$,
since the case $-1 \le x < 0$ is the same.
Let $x=:\cos \theta$ and denote
$I_{[a, b]}:= \int_a^b|L_n^{\a,\b} (x,y)|^p w_{\a,\b}(y) dy$.
We first estimate $I_{[-1/2, 1]}$.
Using Theorem~\ref{thm:3.1} and (\ref{w-sim1})
and substituting $y=\cos\phi$, we get
$$
I_{[-1/2, 1]} \le
c \frac{n^p}{(w_{\a,\b}(n;x))^{p/2}}
\int_0^{2\pi/3} \frac{(\sin \phi)^{2\a+1} d\phi}
{(\sin \phi +n^{-1})^{(\a +1/2)p} (1+n|\theta - \phi|)^{\sigma p}}
$$
with $\sigma= k-2\a-2\b-3$.
Denoting the last integral by $J$, we use the relation
$\sin \phi \sim \phi$ for $0 \le \phi \le 2\pi/3$
and apply the substitution $u=n\theta$ to obtain
\begin{align*}
 J \le & \,c \frac{n^{(\a+1/2) p}}{n^{2\a+1}} \int_0^{2\pi/3} \frac{d\phi}
        {(1+ n \phi)^{(\a+1/2) (p-2)} (1+n|\theta - \phi|)^{\sigma p}}\\
 \le & \,c \frac{n^{(\a+1/2) p}}{n^{2\a +2}} \int_0^{2n\pi/3} \frac{d u}
        {(1+ u)^{(\a+1/2) (p-2)} (1+|u-n \theta|)^{\sigma p} }.
\end{align*}

If $p \ge 2$, then we split the last integral into two integrals: one over
$[0,n \theta/2]$ and the other over $[n \theta/2, 2n\pi/3]$. We have
$$
 \int_0^{n \theta/2} \le \frac{c}{ (1+n\theta)^{\sigma p} }
     \int_0^{n \theta/2} \frac{du} {(1+ u)^{(\a+1/2) (p-2)}}
    \le \frac{c}{(1+n\theta)^{\sigma p -1} }
$$
and
$$
  \int_{n \theta/2}^{2n\pi/3} \le
     \frac{c}{(1+n\theta)^{(\a+1/2)(p-2)}}
      \int_{n\theta /2}^{2n\pi/3} \frac{du}{(1+|u-n \theta|)^{\sigma p}}
   \le \frac{c}{(1+n\theta)^{(\a+1/2)(p-2)}},
$$
using that the last integral is bounded by
$\int_{-\infty}^\infty (1+|u|)^{-\sigma p}du\le c$.

The assumption $k \ge 3\a+3\b+5$ yields $\sigma p -1 \ge \a(p-2)$
and hence
\begin{align*}
 J \le c\, \frac{n^{(\a+1/2)p}}{n^{2\a +2}(1+n\theta)^{(\a+1/2)(p-2)}}
     & \le \frac{c}{n (\theta+n^{-1})^{(\a+1/2)(p-2)}} \\
     & \le \frac{c}{n (w_{\a,\b}(n; \cos\theta))^{p/2-1}}.
\end{align*}

If $p  < 2$, then we have
\begin{align*}
 J \le & \,c \frac{n^{(\a+1/2) p}}{n^{2\a+2}}
    \int_0^{2n \pi/3} \frac{ (1+ u )^{(\a+1/2) (2-p)} d u}
        {(1+|u-n \theta|)^{\sigma p}}\\
   = & \,c \frac{n^{(\a+1/2) p}}{n^{2\a +2}}
    \int_{-n \theta}^{n (2\pi/3-\theta)}
     \frac{(1+ v + n\theta)^{(\a+1/2)(2-p)} d v} {(1+|v|)^{\sigma p} }.
\end{align*}
Using the inequality $(A+B)^a \le 2^a (A^a + B^a)$ we
conclude that
\begin{align*}
 J \le & \,c \frac{n^{(\a+1/2) p}}{n^{2\a+2}}  \left[
     \int_{-\infty}^\infty \frac{d v} {(1+|v|)^{\sigma p -(\a+1/2)(2-p)}}
    + (n \theta)^{(\a+1/2)(2-p)}
     \int_{-\infty}^\infty \frac{d v} {(1+|v|)^{\sigma p}}\right] \\
    \le & \,c \frac{n^{(\a+1/2) p}}{n^{2\a +2}}(1+n\theta)^{(\a+1/2)(2-p)}
    \le \frac{c}{n (w_{\a,\b}(n; \cos\theta))^{p/2-1}},
\end{align*}
which is the same as in the case of $p \ge 2$.
Here we have used that $k \ge 3\a+3\b+5$, which yields
$\sigma p -(\a+1/2)(2-p) > 1$.

Putting these estimates together, we conclude that
\begin{equation}\label{estim-I1}
I_{[-1/2,1]} \le c \,\frac{n^pJ}{(w_{\a,\b}(n;x))^{p/2}}
 \le c \,\frac{n^{p-1}}{(w_{\a,\b}(n;x))^{p-1}}.
\end{equation}

To estimate $I_{[-1, -1/2]}$ we again apply the substitution
$y=\cos \phi$.
Notice that since $0\le \theta\le \pi/2$ and $2\pi/3\le \phi\le \pi$,
we have $|\theta-\phi| \ge \pi/6$.
Employing Theorem~\ref{thm:3.1} and (\ref{w-sim2}), we obtain as above
\begin{align*}
I_{[-1,-1/2]} & \le c \frac{n^p}{(w_{\a,\b}(n;x))^{p/2}}
\frac{ n^{(\b+1/2) p}}{n^{\sigma p}} \int_{2\pi/3}^\pi
\frac{ (\sin \phi)^{2\b+1}} {(1+n \sin \phi)^{(\b+1/2) p}} d \phi \\
& \le c \frac{n^p}{(w_{\a,\b}(n;x))^{p/2}}
\frac{ n^{(\b+1/2) p}}{n^{\sigma p}}
\frac{1}{n^{(\b+1/2)(p-1)}}
 \le c \frac{n^{p-1}}{(w_{\a,\b}(n;x))^{p-1}},
\end{align*}
using the trivial inequalities
$n^{-(2\a+1)}\le w_{\a,\b}(n;\cos\theta)\le c$ for $0\le \theta\le \pi/2$
and the fact that $k\ge 3\a+3\b+5$.
The last estimate coupled with (\ref{estim-I1}) yields (\ref{est-Lp-norm}).
The~proof of Proposition~\ref{prop:3.2} is complete.
\end{proof}

\section{Tight polynomial frames in 
$L^2(w_{\a,\b})$}
\label{frames}
\setcounter{equation}{0}

We now turn to the construction of polynomial frames.
In addition to the localized polynomials from \S\ref{local-polyn},
we shall need the Gaussian quadrature formula. Let $\xi_j = \cos \theta_j$,
$1 \le j \le n$, denote the zeros of the Jacobi polynomial
$P_n^{(\a,\b)}(t)$, ordered so that
$$
0=:\theta_0 < \theta_1 < \dots < \theta_n < \theta_{n+1}:=\pi
$$
Let $\lambda_n(t)$ be the Christoffel function and $b_\nu =
\lambda_n(\xi_\nu)$. It is known that
$$
   \theta_{\nu+1} - \theta_\nu \sim n^{-1}
   ~ \mbox{and hence $~\theta_\nu \sim \nu n^{-1}$ $(1 \le \nu \le n)$;}
$$
and also
$$
    b_\nu \sim  n^{-1} w_{\a, \b}(\xi_\nu)(1-\xi_\nu^2)^{1/2} \sim
     n^{-1} w_{\a,\b}(n; \xi_\nu).
$$
Here the constants of equivalence depend only on $\a,\b$
(cf. \cite[p. 282]{SV}). These quantities appear in the well-known
Gaussian quadrature:

\begin{prop}\label{p:quadrature}
For each $n \ge 1$, the quadrature
\begin{equation}\label{quadrature}
c_{\a,\b}\int_{-1}^1 f(t)w_{\a,\b}(t)dt \sim \sum_{\nu=1}^n b_\nu f(\xi_\nu)
\end{equation}
is exact for all polynomials of degree $2n-1$.
\end{prop}

We shall utilize Proposition \ref{p:quadrature} with $n=2^{j}$. Let us denote
by $\xi_{j,\nu}$ and $b_{j,\nu}$ ($\nu=1, 2,\dots, 2^{j}$) the
knots and the coefficients of the Gaussian quadrature when $n = 2^j$, and
set
$$
\cX_j:=\{\xi_{j,\nu}: \nu=1, 2,\dots, 2^{j}\}, \quad
\hbox{and}\quad b_\xi:=b_{j,\nu} \quad \hbox{if $\xi=\xi_{j,\nu}$}.
$$
Then the Gaussian quadrature
\begin{equation}\label{quadrature1}
c_{\a,\b}\int_{-1}^1 f(t)w_{\a,\b}(t)dt \sim \sum_{\xi \in \cX_j} b_\xi f(\xi)
\end{equation}
is exact for all polynomials of degree at most $2^{j+1} -1$.

Let $\wh a$ satisfy the conditions:
\begin{equation}\label{a1}
\wh a \in C^\infty(\RR), \quad \wh a \ge 0,
 \quad \supp \wh a \subset [1/2, 2],
\end{equation}
\begin{equation}\label{a2}
\wh a(t)>c>0, \qquad \mbox{if $t \in [3/5, 5/3]$},
\end{equation}
\begin{equation}\label{a3}
\wh a^2(t) + \wh a^2(2t) =1,
\qquad \mbox{if $t \in [1/2, 1]$.}
\end{equation}
It follows from these conditions that
\begin{equation}\label{a4}
\sum_{\nu=0}^\infty \wh a^2(2^{-\nu}t)= 1,
\qquad t \in [1, \infty).
\end{equation}

It is easy to construct functions that satisfy properties
\eqref{a1}-\eqref{a3}. Indeed, a standard construction (e.g. in wavelets) shows
that there exists a function $g$ with the properties:
$g \in C^\infty(\RR)$, $\supp g = [-1, 1]$,
$g(t) > 0$ on $(-1, 1)$,
$g(-t) = g(t)$, $g(0) =1$
and $|g(t)|^2 + |g(t+1)|^2 =1$ on  $[-1, 0]$. Then the function
$\wh a(t) := g(\log_2t)$ has the required properties.

We introduce the orthonormal Jacobi polynomials $\PP_n^{(\a,\b)}$, which
are given by
$$
\PP_n^{(\a,\b)} (x) := (h_n^{(\a,\b)})^{-1/2} P_n^{(\a,\b)}(x).
$$
The sequence $\{\PP_\nu^{(\a,\b)}\}_{\nu=0}^\infty$ is an orthonormal basis in
$L^2(w_{\a,\b},[-1, 1])$ with inner product
$$
 \langle f, g\rangle: = c_{\a,\b} \int_{-1}^1 f(t)\overline{g(t)}
   w_{\a,\b}(t) dt.
$$

Assuming that $\wh a$ satisfies conditions (\ref{a1})-(\ref{a3}), we write
$L_0^{\a,\b} := 1$ and
\begin{equation}\label{def.Lambdaj}
L_j^{\a,\b} (x, y) := \sum_{\nu=0}^\infty
\wh a \Big(\frac{\nu}{2^{j-1}}\Big)\PP_\nu^{(\a,\b)}(x)\PP_\nu^{(\a,\b)}(y),
\quad  j=1, 2, \dots.
\end{equation}
We now define our frame elements ({\bf needlets}) by
\begin{equation}\label{def.frame}
\psi_\xi(x):= \sqrt{b_\xi}\cdot L_j^{\a,\b}(x, \xi)
\quad \mbox{for}\quad
\xi\in \cX_j, ~ j=0, 1, \dots,
\end{equation}
where $\cX_j$ is the set of knots and the $b_\xi$'s are the weights
in the Gaussian quadrature~\eqref{quadrature1}.

Denote $\cX:=\cup_{j=0}^\infty \cX_j$, where every two points
$\xi, \eta \in \cX$
from different levels $\cX_j\ne\cX_k$
are considered to be different elements of $\cX$ even if they coincide.
We use $\cX$ as an index set in the definition of the needlet system
$$
\Psi:=\{\psi_\xi\}_{\xi\in\cX},
$$
which as will be shown in Theorem~\ref{t:frame}
below is a tight frame in $L^2(w_{\a,\b})$.

We next show that each $\psi_\xi$ has nearly exponential
(faster than any polynomial) rate of decay away from $\xi$.
This property of needlets is critical for their implementation
to decomposition of spaces other than $L^2(w_{\a,\b})$.


\begin{thm}\label{p:localization}
For any $k \ge 1$ there exists a constant $c_k$ depending only on $k$,
$\a,\b$, and $\wh a$ such that for $\xi \in X_j$, $j=1, 2, \dots$,
\begin{equation}\label{loc-frame}
|\psi_\xi(\cos \theta)|
\le c_k \frac{2^{j/2}}{\sqrt{w_{\a,\b}(2^j;\cos\theta)}
    (1+2^j|\theta-\arccos \xi|)^k},
\qquad \theta\in [0, \pi],
\end{equation}
and, consequently, for $\nu=1, 2, \dots, 2^{j}$,
\begin{equation}\label{loc-frame1}
|\psi_{\xi_{j,\nu}}(\cos \theta)|
\le c_k\frac{2^{j/2}}{\sqrt{w_{\a,\b}(2^j;\cos\theta)}
(1+2^j|\theta-\nu  2^{-j}|)^k},
\qquad \theta\in [0, \pi].
\end{equation}
\end{thm}

\begin{proof}
Estimates (\ref{loc-frame})-(\ref{loc-frame1}) follow immediately
from the estimate in Theorem~\ref{thm:3.1}, taking into account the simple
fact that $1+2^j|\theta\pm c2^{-j}| \sim 1+2^{j}\theta$ if
$0\le \theta\le \pi$ and the fact that $b_\xi \sim 2^{-j}
w_{\a,\b}(2^{-j};\xi)$. Since by assumption $\wh a$ is in $C^\infty$,
estimates (\ref{loc-frame})-(\ref{loc-frame1}) holds for every $k \ge 1$.
\end{proof}


\noindent
{\em Remark.}
It is worthwhile to mention that the estimate (\ref{loc-frame})
or (\ref{loc-frame1}) (also Proposition~\ref{prop:3.2} with $p = 2$)
immediately yields
$\|\psi_{\xi_{j,\nu}}\|_{L^2(w_{\a,\b})}\le c$.


\medskip

Our next theorem shows that $\Psi$ is a tight frame in $L^2(w_{\a,\b})$.


\begin{thm}\label{t:frame}
If $f\in L^2(w_{\a,\b})$, then
\begin{equation}\label{frame1}
f =\sum_{j=0}^\infty\sum_{\nu=1}^{2^j}
  \langle f, \psi_{\xi_{j,\nu}} \rangle \psi_{\xi_{j,\nu}}
= \sum_{\xi\in \cX} \langle f, \psi_\xi\rangle \psi_{\xi}
\quad\mbox{in $L^2(w_{\a,\b})$}
\end{equation}
and
\begin{equation}\label{frame2}
\|f\|_{L^2(w_{\a,\b})}
=\Big(\sum_{\xi\in \cX} |\langle f, \psi_\xi\rangle|^2\Big)^{1/2}.
\end{equation}
\end{thm}

\begin{proof}
Suppose $f\in L^2(w_{\a,\b})$. Then
$f=\sum_{\nu=0}^\infty d_\nu(f)\PP_\nu^{(\a,\b)}$,
where $d_\nu(f):= \langle f, \PP_\nu^{(\a,\b)} \rangle$.
Here and elsewhere in the following the convergence is in $L^2(w_{\a,\b})$.
Denote briefly
$$
(L_j*f)(x):= c_{\a,\b} \int_{-1}^1 L_j^{\a,\b}(x, t)f(t)w_{\a,\b}(t)dt,
$$
where $L_j^{\a,\b}$ are defined in (\ref{def.Lambdaj}).
One can regard this as a nonstandard convolution which is evidently
associative but not commutative.
We next show that the following Calderon type decomposition holds
\begin{equation}\label{L.repr}
f=\sum_{j=0}^\infty L_j^{\a,\b}*L_j^{\a,\b}*f.
\end{equation}
Indeed, clearly for $j\ge 1$,
$$
L_j^{\a,\b}*f = \sum_{\nu=1}^{2^j} \wh a \Big(\frac{\nu}{2^{j-1}}\Big)
d_\nu(f)\PP_\nu^{(\a,\b)}
$$
and then
$$
L_j^{\a,\b}*L_j^{\a,\b}*f
= \sum_{\nu=1}^{2^j} \wh a^2 \Big(\frac{\nu}{2^{j-1}}\Big)
d_\nu(f)\PP_\nu^{(\a,\b)}.
$$
Consequently,
\begin{align*}
\sum_{j=0}^\infty L_j^{\a,\b}*L_j^{\a,\b}*f
=& \,d_0(f)\PP_0^{(\a,\b)} +\sum_{j=1}^\infty\sum_{\nu=1}^{2^j}
\wh a^2 \Big(\frac{\nu}{2^{j-1}}\Big)d_\nu(f)\PP_\nu^{(\a,\b)}\\
=& \,d_0(f)\PP_0^{(\a,\b)} +\sum_{\nu=1}^\infty\sum_{j=1}^\infty
\wh a^2 \Big(\frac{\nu}{2^{j-1}}\Big)d_\nu(f)\PP_\nu^{(\a,\b)}\\
=& \sum_{\nu=0}^\infty d_\nu(f)\PP_\nu^{(\a,\b)} = f,
\end{align*}
where we used (\ref{a4}). Thus (\ref{L.repr}) is established.

As we already mentioned
$L_j^{\a,\b}*L_j^{\a,\b}*f= (L_j^{\a,\b}*L_j^{\a,\b})*f$.
Since $L_j^{\a,\b}(x, t)L_j^{\a,\b}(t, y)$
is a polynomial of degree at most $2^{j+1}-1$ in $t$
and the quadrature formula (\ref{quadrature1})
is exact for such polynomials, we have
\begin{align*}
(L_j^{\a,\b}*L_j^{\a,\b})(x, y)
=& \, c_{\a,\b} \int_{-1}^1 L_j^{\a,\b}(x, t)L_j^{\a,\b}(t, y)w_{\a,\b}(t)dt\\
=& \, \sum_{\xi\in \cX_j} \sqrt{b_\xi}\cdot L_j^{\a,\b}(x, \xi)
\sqrt{b_\xi}\cdot L_j^{\a,\b}(\xi, y)\\
=& \, \sum_{\xi\in \cX_j} \psi_\xi(x) \psi_\xi(y).
\end{align*}
Here we also used that $L_j^{\a,\b}(x, y)=L_j^{\a,\b}(y, x)$ and
the definition of the frame elements in (\ref{def.frame}).
Consequently,
$$
L_j^{\a,\b}*L_j^{\a,\b}*f
= \sum_{\xi\in \cX_j} \langle f, \psi_\xi \rangle \psi_\xi,
$$
which along with (\ref{L.repr}) yields (\ref{frame1}).

For the proof of (\ref{frame2}) we denote
$
F_Jf:= \sum_{j=0}^J\sum_{\xi\in \cX_j} \langle f, \psi_\xi \rangle \psi_\xi.
$
Evidently,
$$
\langle f, F_Jf \rangle
= \sum_{j=0}^J\sum_{\xi\in \cX_j} |\langle f, \psi_\xi \rangle|^2
$$
and passing to the limit as $J\to \infty$ we arrive at (\ref{frame2}).
\end{proof}

\end{document}